\documentclass[11pt,twoside,final]{amsart}

\newtheorem{definition}{Definition}[section]
\newtheorem{theorem}{Theorem}[section]
\newtheorem{proposition}[theorem]{Proposition}
\newtheorem{lemma}[theorem]{Lemma}
\newtheorem{corollary}[theorem]{Corollary}
\newtheorem{remark}[theorem]{Remark}

\usepackage[english]{babel}
\usepackage{graphicx,epstopdf,epsfig}
\usepackage{amsfonts,epsfig,fancyhdr,graphics,amsmath,amssymb}

\title{Bireflectionality in special orthogonal groups}

\newcommand {\ch }{\operatorname{char}}
\newcommand {\SOG }{\operatorname{SO}}
\newcommand {\OG }{\mathrm{O}}
\newcommand {\GL }{\mathrm{GL}}
\newcommand {\Bahn }{\mathrm{B}}
\newcommand {\Fix }{\mathrm{Fix}}
\newcommand {\GF }{\mathrm{GF}}
\newcommand {\SpG }{\mathrm{Sp}}
\newcommand {\Idm }{\mathrm{I}}
\newcommand {\Cent }{\mathrm{Cent}}
\newcommand {\Arf }{\operatorname{Arf}}
\newcommand {\rad }{\operatorname{rad}}
\newcommand {\Wind }{\operatorname{ind}}

\usepackage{ifdraft}
\ifdraft{\today}{}
\ifdraft{\usepackage[outer]{showlabels}
	\usepackage{todonotes}}{}
\usepackage[pagebackref,colorlinks,citecolor=red,urlcolor=blue,linkcolor=green, bookmarks=false,hypertexnames=true]{hyperref}
\usepackage{orcidlink}
\date{September, 28, 2024}

\begin{document}

\bibliographystyle{plain}

\setcounter{page}{1}

\thispagestyle{empty}


\author{
Klaus Nielsen}\, \orcidlink{0009-0002-7676-2944}
\email{klaus@nielsen-kiel.de}
\maketitle

\begin{abstract}
It is shown that a transformation in the special orthogonal group $\SOG(V,Q)$ of a nondefective quadratic space $(V,Q)$ over a field of arbitrary characteristic  is bireflectional (product of 2 involutions) if and only if it is reversible (conjugate to its inverse). Furthermore, all elements of $\SOG(V,Q)$ are bireflectional if and only if $\dim V \not \equiv 2 \mod 4$ or $(V,Q)$ is a hyperbolic plane over $\GF(2)$ or $\GF(3)$.
\end{abstract}

\section{Introduction and main results} \label{intro-sec}

Following H. Wiener\footnote{H. Wiener, Über die aus zwei Spiegelungen zusammengesetzten Verwandtschaften.
	Ber. Verh. kgl. Sächs. Ges. Wiss. Leipzig. Math.-phys. Cl. 43: 644-673 (1891)}, we call an element of a group $G$ bireflectional if it is a product of 2 involutions of $G$.
We say that an element of $G$ is reversible if it is conjugate to its inverse. Clearly, a bireflectional element is reversible.

	A group  is bireflectional if all its elements are bireflectional.
 Bireflectional elements and groups are also called strongly real by some authors;
 see e.g. \cite{KTV-2023} and \cite{Ramo-2011}.
 
We are interested in classifying bireflectional elements in some subgroups of orthogonal groups. 
In \cite{KN-2010}, F. Knüppel and the author considered bireflectionality in the commutator subgroup of an orthogonal group over the reals. In a forthcoming paper
we classify bireflectional and reversible elements in the commutator subgroup of an orthogonal group over a finite field.

In this paper, we study bireflectionality in the special orthogonal group. 
Let $(V,Q)$ be a quadratic space, where $V$ is a vector space of finite dimension over a field $K$ and $Q$ is a nondefective quadratic form on $V$.
The special orthogonal group $\SOG(V) = \SOG(V,Q)$ is the subgroup of all orthogonal transformations $\varphi$ of $(V, Q)$ with $\dim V(\varphi -1) \equiv 0 \mod 2$.  

It is well-known that the full orthogonal group $\OG(V) = \OG(V,Q)$ is bireflectional by theorems of Wonenburger \cite{Wonenburger-1966} and Ellers \& Nolte \cite{EllersNolte-1982} and also by Gow \cite{Gow-1981}. 
In \cite{KN-1987}, F. Knüppel and the author proved the following
theorem.

\begin{theorem}(\cite[Theorem A]{KN-1987}\cite[Theorem B]{KN-1987})\label{old-theorem}
	Let $\ch K \ne 2$.
	\begin{enumerate}
	 \item The special orthogonal group $\SOG(V)$ is bireflectional if and only if $\dim V \not \equiv 2 \mod 4$ or $(V,Q)$ is a hyperbolic space over  $\GF(3)$.
	 \item Let $\dim V \equiv 2 \mod 4$ and $\varphi \in \SOG(V)$. Then $\varphi$ is bireflectional if and only if
	 $\varphi$ has an orthogonal summand of odd dimension.
\end{enumerate}
\end{theorem}

 Only after the publication of \cite{KN-1987}, we considered the  case $\ch K = 2$ and obtained the following theorem.

\begin{theorem} \label{theorem-1}
	Let $\ch K = 2$.
	\begin{enumerate}
		\item Let $\dim V \ge 4$. The special orthogonal group $\SOG(V)$ is bireflectional if and only if $\dim V  \equiv 0 \mod 4$.
		\item Let $\dim V \equiv 2 \mod 4$ and $\varphi \in \SOG(V)$. Then $\varphi$ is bireflectional if and only if
		$\varphi$ has a unipotent orthogonal summand $\psi$, where 
		\begin{enumerate}
			\item $\psi$ is  cyclic, or
			\item $\psi$ is orthogonally indecomposable and  $\dim \psi \equiv 2 \mod 4$.
		\end{enumerate}
	\end{enumerate}
\end{theorem}
We have never published it so far. And to my knowlodge, nobody before has considered the $\ch K = 2$ case. Except in the case that $K$ is finite. 
J. Rämö \cite{Ramo-2011} proved that special orthogonal group $\SOG(V)$ over a finite field of characteristic 2  is bireflectional if and only if $\dim V \equiv 0 \mod 4$.
Furthermore, she showed that all unipotent elements of $\SOG(V)$ are bireflectional. In a recent paper, Kim et al. \cite[Theorem 12]{KTV-2023} showed that all reversibles of $\SOG(V)$ are already bireflectional. They also gave a characterization of bireflectional elements. But their result seems to be not quite correct.

Here we give a new proof of theorem \ref{theorem-1} which is slightly different from our old one.
As in the old proof, we need some elementary facts
about orthogonally indecomposable orthogonal transformations. These can be found in papers of Zassenhaus \cite{Zassenhaus1958b}, Wonenburger \cite{Wonenburger-1966}, and Huppert \cite{Huppert-1980a}. Or more general in \cite{GKN-2008}.

An orthogonal transformation $\varphi \in \OG(V)$ is  orthogonally indecomposable if $V$ has no proper nondefective $\varphi$-invariant subspaces. 
According to Huppert \cite[1.7 Satz]{Huppert-1980a}, an orthogonally indecomposable transformation is either\footnote{Our type enumeration follows Huppert \cite{Huppert-1980a}. In \cite{Huppert-1990} Huppert uses  a different  enumeration.}
\begin{enumerate}
	\item bicyclic with elementary divisors $e_1 = e_2 =(x \pm 1)^m$
	(type 1), or
	\item indecomposable (type 2), or
	\item cyclic and biprimary ( type 3).
\end{enumerate}

\begin{definition}
	Let $\varphi \in \OG(V)$. We say that $\varphi$ is hyperbolic if $V = U \oplus W$ for some totally isotropic $\varphi$-invariant subspaces $U$ and $W$ of $V$.
\end{definition}

Clearly, if $\OG(V)$ contains hyperbolic transformations, then $V$ must be hyperbolic. In \cite[2.3 Satz]{Huppert-1980a}, \cite[2.4 Satz]{Huppert-1980a}, and \cite[3.2 Satz]{Huppert-1980b} Huppert proved the following.

\begin{proposition} \label{prop-1}
	Let $\ch K \ne 2$.
	Let $\varphi \in \OG(V)$ be orthogonally indecomposable
	of type 1. Then 
	\begin{enumerate}
		\item If $\varphi$ is of type 1, then $\dim V \equiv 0 \mod 4$ and $\varphi$ is hyperbolic.
		\item $\dim V \equiv 1 \mod 2$ if and only if $\varphi$  or $-\varphi$ is unipotent
	\end{enumerate}
\end{proposition}

\begin{corollary}
	Let $\ch K \ne 2$.
	Let $\varphi \in \OG(V)$ be orthogonally indecomposable
	of type 1. There exists a cyclic transformation $\psi \in \OG(V)$  with minimal polynomial $(x^2+1)^m$. Further, 
	\begin{enumerate}
		\item If $\varphi$ and $\zeta \in \OG(V)$ are similar, then they are conjugate.
		\item If $-\varphi$ is unipotent, then $\varphi \sim \psi^2$. 
		\item If $\varphi$ is unipotent, then $\varphi  \sim \psi_U$, the unipotent factor in the multiplicative Jordan-Chevalley decomposition of $\psi$.
		\end{enumerate}
\end{corollary}

The existence of $\psi$ is shown in Huppert \cite[3.2 Satz]{Huppert-1980b}.
Theorem \ref{prop-1} does not hold if $\ch K = 2$.
Huppert \cite{Huppert-1980a}, 2.5 Hilfssatz] shows that there exists a
nonhyperbolic  orthogonally indecomposable transformation of type 1 of dimension 6. We give further examples.

\begin{remark} \label{remark-1}
	Let $(V, h)$ be a nondegenerate unitary vector space over a field of characteristic 2. 
	Let $\varphi$ be a unitary transformation of $(V, h)$. Then $\varphi \in \OG(V,Q_h)$, where $Q_h$ is the trace form of $h$. If $\dim V$ is odd, and   $\varphi$ is  cyclic and unipotent, then $\dim (V,Q_h) \equiv 2 \mod 4$, and $\varphi$ is orthogonally indecomposable of type 1 w.r.t. $Q_h$. By \cite[Remark 10.1.4]{WScharlau}, 
	$(V,Q_h)$ has a nontrivial Arf invariant. Hence $(V,Q_h)$
	is not hyperbolic. This follows also from Bünger \cite[Bemerkung 2.2.4]{Bunger-1997}.
\end{remark}

\begin{remark} \label{remark-1a}
	Let $(V,Q)$ be a nondefective quadratic space over a field of characteristic 2. Assume that $V$ has dimension $4m+2$ and Witt index $2m$.
	Let $\psi \in \OG(V)$ be cyclic.
	\begin{enumerate}
		\item If $\psi \in \OG(V)$ is unipotent, then  $\psi^2$ is orthogonally indecomposable of type 1 and not hyperbolic.
		\item If $\psi$ has minimal polynomial $p^{2m+1}$, where $p$ is irreducible of degree 2, then  the unipotent factor in the multiplicative Jordan-Chevalley decomposition of $\psi$ is orthogonally indecomposable of type 1 and not hyperbolic.
	\end{enumerate}
\end{remark}

\begin{definition}
	{\rm Let $\ch K = 2$. An  orthogonally indecomposable orthogonal transformation of type 1 is of  {\em type 1o}, resp. {\em type 1e} if its minimal polynomial has odd (even) degree.}
\end{definition}

Our new proof of theorem \ref{theorem-1}  is inspired by a paper of Gonshaw et al.\footnote{There is a  minor typo on p. 95, corrected by Giovanni De Franceschi, Centralizers and conjugacy classes in finite classical groups. PhD
	Thesis, University of Auckland (2018)} \cite{Gonshaw_Liebeck_OBrien}
The authors consider unipotent conjugacy classes in finite classical groups. It follows from their paper that  the conjugacy class of a unipotent ortho\-gonally indecomposable  transformation $\varphi \in \OG(V,Q,K)$ is determined by its Arf invariant (i.e the Arf invariant of $(V,Q)$) if $K$ is a finite field of characteristic 2. Consequently,  $\varphi$ has a unipotent cyclic square root $\Psi$ if $\varphi$ is of type 1o. In this case, $\psi = \rho \sigma$ is the product of 2 orthogonal involutions, and $\varphi = \rho \rho^{\sigma} = \sigma^{\rho} \sigma$. Now either $\rho \not \in \SOG(V)$ or $\sigma \not \in \SOG(V)$ as $\psi \not \in \SOG(V)$. 

What can be said in the case of infinite fields? In general, 
the conjugacy class of a cyclic unipotent transformation additionly depends  on its spinor norm (which is trivial if $K$ is perfect).
[If $\ch K \ne 2$ a cyclic unipotent  conjugacy class $\Omega$ is determined by the discriminant of $(V,Q)$, which is the spinor norm of $-\Omega$, alone.]. 
Bünger \cite[Satz.15.5]{Bunger-1997} proved that $\varphi$ is always hyperbolic if $\varphi$ is orthogonally indecomposable of type 1e.

We show

\begin{proposition}\label{prop-2}
	Let $\ch K = 2$.
	Let $\varphi \in \OG(V,Q)$ be orthogonally indecomposable
	of type 1. 
	\begin{enumerate}
		\item Let $\varphi$ be of type 1e. Then
		\begin{enumerate}
			\item $\varphi$ is hyperbolic. 
		\item $\varphi$ is the unipotent factor in the multiplicative Jordan-Chevalley decomposition of a cyclic  transformation in $\OG(V,Q)$.
	\end{enumerate}
		\item $\varphi$ has a square root in $\OG(V,Q)$ if and only if $\varphi$ is of type 1o.
	\end{enumerate}
\end{proposition}

It follows that the spinor norm of an orthogonally indecomposable transformation of type 1o (and 1e) is trivial. But we don't know whether its conjugacy class is already determined by its Arf invariant.

Finally, we also want to prove:

\begin{theorem} \label{theorem-2}
	All reversible elements of $\SOG(V,Q)$ are bireflectional.
\end{theorem}
It is shown in \cite[Theorem 4.1]{KN-2010} that this true if $\ch K \ne 2$.

\section{Notations}
Let $V$ be a vector space of finite dimension $n$ over a field $K$.
For a linear mapping $\varphi$ of $V$ let $\Bahn^j(\varphi)$ denote the image and $\Fix^j(\varphi)$ the kernel of $(\varphi -1)^j$. Let $\Fix^{\infty}(\varphi) :=  \Fix^n(\varphi)$ denote the Fitting null space of $\varphi -1$ and $\Bahn^{\infty}(\varphi): = 
\Bahn^n(\varphi)$ the Fitting one space of $(\varphi -1)^j$. The space $\Bahn(\varphi) := \Bahn^1(\varphi)$ is the path or residual space of $\varphi$.
And $\Fix(\varphi) := \Fix^1(\varphi)$ is the fix space of $\varphi$. By $\mu(\varphi)$ we denote the minimal polynomial of 
$\varphi$.
Clearly, the spaces $\Bahn^j(\varphi)$ and $\Fix^j(\varphi)$ are $\varphi$-invariant. The centralizer of $\varphi$ is denoted by $\Cent(\varphi)$.

Let $(V,Q)$ be a quadratic space, i.e $Q$ is a quadratic form on $V$. Let $f_Q$ denote its associated bilinear form defined by $f_Q(u,w) =Q(u+w) - Q(u) - Q(w)$.
We always assume that $Q$ is nondefective; i.e $f_Q$ is nondegenerate.

 A quadratic subspace space $(W,Q)$ of $(V,Q)$ is totally isotropic if $Q(w) = 0$ for all $w \in W$. It is totally degenerate if $W \le  W^{\perp}$. 
 If $\rad W = W \cap W^{\perp}$ is totally isotropic, then $Q$ induces a nondefective quadratic form on $W/\rad W$.

 If $Q$ is hyperbolic, i.e.  $Q$ has Witt index $\Wind Q = \frac{\dim V}{2}$, then the Arf invariant $\Arf Q$ of $Q$ is zero.
For basic facts about quadratic forms we refer to the books of Scharlau \cite{WScharlau} and Kaplansky \cite{Kaplansky-1974}. We collect some well-known facts.

\begin{remark} \label{remark-2}
	Let $\varphi \in \OG(V)$. Then 
	\begin{enumerate}
		\item If $W$ is a subspace of $V$, then $\rad W \cap \Bahn(\varphi)$ is totally isotropic.
		\item $\Bahn^j(\varphi)^{\perp} = \Fix^j(\varphi)$.
		\item $V = \Bahn^{\infty}(\varphi) \oplus \Fix^{\infty}(\varphi)$.
		\item If $\Bahn(\varphi)$ is totally isotropic, then $\dim \Bahn(\varphi)$ is even.
	\end{enumerate}	
\end{remark}

If $b \in V$ is anisotropic, then there exists a unique transformation $\sigma_b \in \OG(V,Q)$ with $\Bahn(\sigma) = \langle b\rangle$, the reflection along $b$.
We have 
\[
v \sigma_b = v - \frac{f_Q(v,b)}{Q(b)} b.
\]
The following is well-known and shows that $\SOG(V)$ is actually a group and coincides with the usual special orthogonal group if $\ch K \ne 2$. It is  also used to prove the theorem of Cartan-Dieudonn\'e-Scherk. For the convenience of the reader, we  provide a short proof.
\begin{remark} \label{remark-3}
	Let  $\varphi \in \OG(V)$, and let $\sigma \in \OG(V,Q)$ be a reflection. Then either 
	\begin{enumerate}
	\item $\Bahn(\varphi \sigma) = \Bahn(\varphi) \oplus \Bahn(\sigma)$, or 
	\item $\Bahn(\varphi) = \Bahn(\varphi \sigma) \oplus \Bahn(\sigma)$.
\end{enumerate}	
\end{remark}

\begin{proof}
	1: Assume that $\Bahn(\sigma) \cap \Bahn(\varphi) = 0$. Then 
	$\Fix(\varphi \sigma) = \Fix(\varphi ) \cap \Fix(\sigma)$:
	If $u \in \Fix(\varphi \sigma) - \Fix(\varphi)$, then $0 \ne u - u \varphi 
 = u \varphi \sigma	- u \varphi \in \Bahn(\sigma)$, a contradiction.
 By \ref{remark-2}, $\Bahn(\varphi \sigma) = \Bahn(\varphi) \oplus \Bahn(\sigma)$. 
	
	2: 
	Assume that $\Bahn(\sigma) \le \Bahn(\varphi)$. Then $\Fix(\varphi ) \le \Fix(\sigma)$. Hence $\Fix(\varphi ) \le \Fix(\varphi \sigma)$ so that
	$\Bahn(\varphi \sigma) \le \Bahn(\varphi)$. Then  $\Bahn(\varphi \sigma) + \Bahn(\sigma) \le \Bahn(\varphi) = \Bahn(\varphi \sigma \sigma) \le \Bahn(\varphi) + \Bahn(\sigma)$.
	Let $\Bahn(\sigma) = \langle v - v \varphi \rangle$. Then $v \in \Fix(\varphi \sigma)$. Hence   $\Fix(\varphi \sigma) \ne \Fix(\varphi)$ and $\Bahn(\varphi \sigma) \ne \Bahn(\varphi)$.
\end{proof}
\section{Proof of proposition \ref{prop-2}}

Throughout this section we always assume that $\ch K = 2$.

\begin{lemma}                                                 \label{lemma-1}
	Let $\varphi \in \OG(V)$ be cyclic with minimal polynomial
	$(x+1)^{2m}$. Then
	\begin{enumerate}
		\item $\Bahn^{m+1}(\varphi) = \Fix^{m-1}(\varphi)$ is  totally isotropic.
		\item $\Bahn^m(\varphi) = \Fix^m(\varphi)$ is not totally
		isotropic.                   
	\end{enumerate} 
\end{lemma}

\begin{proof}
	1:  $\Bahn^{m}(\varphi)$ is totally degenerate, hence $\Bahn^{m+1}(\varphi) = \Bahn^{m}(\varphi)(\varphi+1)$ is totally isotropic. We prove
	2: Let $Q$ be nondefective. If $m=1$, then $Q(x-x\varphi) = f_Q(x,x\varphi) \ne 0$.
	If $m \geq 2$, then by induction,  
	$\Bahn^{m+1}(\varphi)/\Fix(\varphi)$ is totally isotropic and
	$\Bahn^m(\varphi)/\Fix(\varphi)$ is not. The assertion follows.
\end{proof}

\begin{corollary}                                                 \label{cor-2}
	Let $\varphi \in \OG(V)$ be orthogonally indecomposable. 
	 \begin{enumerate}
	 	\item $V$ is hyperbolic if $\varphi$ is of type 1e
	 	\item $\Bahn(\varphi)/\Fix(\varphi)$ is  orthogonally indecomposable.
	 \end{enumerate}
\end{corollary}

\begin{proof}
	Clearly, $\Bahn(\varphi)/\Fix(\varphi)$ is  orthogonally indecomposable if $\varphi$ is of type 1o. 
So let   $\mu(\varphi) = (x+1)^{2m}$. Then $\Bahn^m(\varphi) = \Fix^m(\varphi)$ is totally degenerate. 
Let $V = U \oplus W$, where $U$ and $W$ are $\varphi$-invariant. 
Then $\Bahn^m(\varphi|{_U})$  and  $\Bahn^m(\varphi|{_W})$ are 
  totally isotropic: If $U$ is totally degenerate, then $\Bahn(\varphi|{_U})$ is totally isotropic. If $U$ not  totally degenerate, then $\rad U$ is totally isotropic, and by \ref{lemma-1}, $\Bahn^{m-r+1}(\varphi|{_{\widetilde{U}}})$ is totally isotropic, where $\widetilde{U} = U/{\rad U}$ and $r = \frac{\dim \rad U}{2}$. Hence 
  $\Bahn^m(\varphi|{_U})$ is totally isotropic.
   Similarly, $\Bahn^m(\varphi|{_W})$ is totally isotropic. Hence 
 $\Bahn^m(\varphi)$ is totally isotropic. By \ref{lemma-1}(2), $\Bahn(\varphi)/\Fix(\varphi)$ must be   orthogonally indecomposable.
\end{proof}

\begin{lemma}                                               \label{lemma-2}
	Let  $\varphi \in \OG(V)$ be orthogonally
	indecomposable of type 1. Assume that  $V$ is hyperbolic.
	Then $\varphi$ is hyperbolic.
\end{lemma}

\begin{proof}
	There is nothing to prove if $\dim V = 2$.
	Now let $\dim V \geq 4$. Then by induction
	\[
	\Bahn(\varphi)/\Fix(\varphi) = [ \langle u \rangle _\varphi + \Fix(\varphi) ]/\Fix(\varphi)
	\oplus [ \langle w \rangle _\varphi + \Fix(\varphi) ]/\Fix(\varphi),
	\]
	where $\langle u \rangle _\varphi$ and  $\langle w \rangle _\varphi$ are
	totally isotropic.  Then
	$\Bahn(\varphi) = \langle u \rangle _\varphi \oplus \langle w \rangle _\varphi$.
	Since $\Fix(\varphi)$ is totally isotropic by \ref{lemma-1}, there exist isotropic
	vectors $x,z \in V$ such that $u = x + x\varphi, w = z + z\varphi$.
	Clearly, $V= \langle x \rangle _\varphi \oplus \langle z \rangle _\varphi$,
	and $\langle x \rangle _\varphi$ and $\langle z \rangle _\varphi$ are
	totally isotropic.
\end{proof}

\begin{corollary}      \label{lemma-3}
	Let  $\varphi \in \OG(V)$ be  orthogonally indecomposable  of type 1. Then $\varphi$ is hyperbolic if one of the following holds
	\begin{itemize}
		\item $K$ is quadratically closed;
		\item  $\varphi$ is of type 1e;
		\item $\dim \rad U \ge 3$ for some $\varphi$-cyclic subspace $U$
		with $\dim U = \frac{\dim V}{2}$.
	\end{itemize}
\end{corollary}

\begin{proof}
	Clearly, $V$ is hyperbolic if $K$ is quadratically closed.
	If $\varphi$ is of type 1e, then $V$ is hyperbolic by \ref{cor-2}. 
	Finally, let $\dim V =4m +2$ and $\dim \rad U \ge 3$. Then $\Bahn^m(\varphi|_U)$ and $\Fix^m(\varphi)$ are totally isotropic by
	\ref{lemma-1}. Hence $\Bahn^m(\varphi|_U) + \Fix^m(\varphi)$ is totally isotropic. So $V$ is hyperbolic.
\end{proof}

\begin{lemma} \label{lemma-4}
	Let $\varphi \in \OG(V)$ be orthogonally indecomposable of type 1o.
	There exist vectors $x, z \in V$ such that
	\begin{enumerate}
		\item $V = \langle x \rangle _{\varphi} \oplus \langle z \rangle _{\varphi}$,
		\item $f_Q (x,x \varphi ^k) = f_Q (z, z \varphi ^k)$ for all $k \in \mathbb{N}$,
		\item $f_Q (x\varphi ^k, z) = f_Q (x, z \varphi ^{k-1})
		$ for all $k \in \mathbb{N}$,
		\item $Q(x) = Q(z)$.
	\end{enumerate}
\end{lemma}

\begin{proof}
	If $\dim V = 2$, then $\varphi = 1$, and the assertion is obviously true.
	Let $\dim V \geq  6$.  By induction,
	\[
	\Bahn(\varphi)/\Fix(\varphi) = [ \langle u \rangle _\varphi + \Fix(\varphi) ]/\Fix(\varphi)
	\oplus [ \langle w \rangle _\varphi + \Fix(\varphi) ]/\Fix(\varphi)
	\]
	for some $u, w \in \Bahn(\varphi)$ with 
	\begin{align*}
		f_Q (u, u\varphi^j) & =  f_Q (w,w\varphi^j), \\
		f_Q (u\varphi ^j, w) & = f_Q (u, w \varphi ^{j-1}), \\
		Q(u) & = Q(w).
	\end{align*}
	 Then $\Bahn(\varphi) = \langle u \rangle _\varphi \oplus
	\langle w \rangle _\varphi$. There exist vectors $x, z \in V$ with
	$u = x(\varphi + 1), w = z(\varphi + 1)$. Clearly,
	$V = \langle x \rangle _\varphi \oplus \langle z \rangle _\varphi$. 
	The equations above  still hold if we replace $w$ by a vector $\widetilde{w}  \in \Fix(\varphi) \oplus \langle w \rangle$. Hence we may replace $z$ by a vector $\widetilde{z} \in \Fix^2(\varphi) \oplus \langle z \rangle$ to achieve that $z \in u^{\perp}$.
	
	 We prove 2:
	The first equation yields
	\[
	f_Q (x, x\varphi^{j-1}) + f_Q (x, x\varphi^{j+1})=f_Q (z,z\varphi^{j-1}) + f_Q (z, z\varphi^{j+1}).
	\]    
In particular for $j = 1$, we have $f_Q (x,x\varphi^2) = f_Q (z,z\varphi^2)$.
From $Q(u) = Q(w)$
we deduce that $f_Q (x,x\varphi) = f_Q (z,z\varphi)$. By induction,
$f_Q (x,x\varphi^k) + f_Q (z,z\varphi^k) = f_Q (x,x\varphi^{k-2}) + f_Q (z,z\varphi^{k-2}) = 0$.

3: 
 Obviously, 3 holds for $k=1$ as $z \in u^{\perp}$.
We have 
\[
f_Q (x\varphi^{j+1}, z) + f_Q (x\varphi^{j-1}, z)= f_Q (x,z\varphi^{j-2}) + f_Q (x, z\varphi^{j}).
\]
In particular for $j=1$ the 2nd equation yields
\[
f_Q (x\varphi^2, z)+ f_Q (x,z) = f_Q (x,z\varphi^{-1})+ f_Q (x,z\varphi) = f_Q (x,z)+ f_Q (x,z\varphi),
 \]
 as $z \in u^{\perp}$. Hence
$f_Q (x\varphi^2, z) = f_Q (x,z\varphi)$.
So 3 also holds for  $k=2$. Now by induction,
\[
f_Q (x\varphi^{k+1}, z) + f_Q (x, z\varphi^{k}) = f_Q (x\varphi^{k-1}, z) + f_Q (x,z\varphi^{k-2}) = 0.
\]
4: Finally, we may replace $z$ with a vector $\widehat{z} \in \Fix(\varphi) \oplus \langle z \rangle$ such that $Q(\widehat{z}) = Q(x)$, and $\widehat{z}$ still satisfies 1-3.
\end{proof}

It follows from 1-3 that $\varphi$ is inverted by an involution $\sigma$ with $\dim \Bahn(\sigma) = \frac{\dim V}{2}$:  
 Define $\sigma \in \GL(V)$ by
$(x \varphi^k) \sigma = z \varphi^{-k}$ and 
		 $(z \varphi^{-k}) \sigma = x \varphi^k$. Then $\sigma$ is involutory, $\sigma \in \OG(V,Q)$, and $\sigma$ reverses $\varphi$.
We used this  in our old proof of  \ref{theorem-1}. The involutions constructed by Ellers and Nolte satisfy $\dim \Bahn(\sigma) = \frac{\dim V}{2} -1$.

\begin{remark} \label{remark-1b}
	Let  $\varphi \in \OG(V)$ be orthogonally
	indecomposable of type 1o. Let $\sigma \in \OG(V)$ be an involution inverting $\varphi$. Then   $V = U \oplus W$ is the direct sum of 2
	$\varphi$- and $\sigma$-invariant subspaces $U$ and $W$ if and only if $\sigma \in \SOG(V)$.
\end{remark}

\begin{proof}
	Let $\dim V = 4m+2$, and let $\tau$ be the involution  $\varphi \sigma$. First let $\sigma \in \SOG(V)$.  Then $\dim \Bahn(\sigma) = \dim \Bahn(\tau) = 2m$. Hence $\Bahn(\varphi) = \Bahn(\sigma) \oplus  \Bahn(\tau)$ so that  $V = \Fix(\sigma) + \Fix(\tau)$.
	So there exist a vector $u \in [\Fix(\sigma) \cup \Fix(\tau)] - \Bahn(\varphi)$ and a vector $w \in [\Fix(\sigma) \cup \Fix(\tau)] - [u(\varphi - 1)^m]^{\perp}$. Put $U = \langle u \rangle_{\varphi}, 
	W = \langle w \rangle_{\varphi}$. Conversely, $\dim \Bahn(\sigma|_U), \dim \Bahn(\sigma|_W) = m$ so that $\dim \Bahn(\sigma) = 2m$.
\end{proof}

\begin{corollary}  \label{cor-3}
	Let $\varphi \in \OG(V)$ be orthogonally indecomposable of type 1. Then $\varphi$ has a square root if and only if $\varphi$ is of type 1o.
\end{corollary}

\begin{proof}
	Let $\varphi$ be of type 1o, and let $x, z$ as in \ref{lemma-4}.
	Define $\psi \in \GL(V)$ by $z\varphi^j \psi = x\varphi^j$ and $x \varphi^j \psi = z \varphi^{j+1}$.
	Then $\psi^2 = \varphi$, and by \ref{lemma-4},
	$\psi \in \OG(V,Q)$. Conversely, let $\psi$ be a square root of $\varphi$. Then $\psi$ must be unipotent and cyclic. 
	Suppose that $\mu(\psi) = (x+1)^{4t}$. By \ref{lemma-1}, $\Fix^t(\varphi) = \Fix^{2t}(\psi)$ is not totally isotropic, in contradiction to \ref{cor-2}.
\end{proof}

\begin{lemma}[Huppert \cite{Huppert-1980a}, 2.2 Hilfssatz ] \label{lemma-5}
	Let $\varphi \in \SpG(V,f_Q)$ be bicyclic with minimal polynomial $(x+1)^t$. Let $V = U \oplus W$, where $U$ and $W$ are degenerate and $\varphi$-invariant.
	Then $\varphi$ is orthogonally indecomposable.
\end{lemma}

	\begin{proof}
		Let $u \in U, w \in W$. Clearly, $(\varphi-1)^{t-1}$ is selfadjoint. Hence 
		$f_Q(u+w, (u+w)(\varphi+1)^{t-1}) =
		f_Q(u, w(\varphi+1)^{t-1}) + f_Q(w, u(\varphi+1)^{t-1}) = 0$.
	\end{proof}
	
	\begin{lemma}                                           \label{lemma-6}
		Let $\varphi, \psi \in \OG(V,Q)$ be orthogonally indecomposable of type 1e. Then  $\varphi$ and $\psi$ are conjugate.
	\end{lemma}
	
	\begin{proof}
		
		Let $V = U \oplus W = X \oplus Z$, where $U, W, X$ and $Z$ are totally isotropic, $U = U\varphi, W = W \varphi$, $X = X\psi$, and $Z = Z \psi$.
		There exists $\alpha \in \OG(V,f)$ with $U \alpha = X$ and $W \alpha = Z$. Then in a suitable basis,
		\[
		\varphi^{\alpha} = \left (\begin{array} {cc} A & 0\\ 0 & A^+ \end{array} \right ),
		\psi = \left (\begin{array} {cc} D & 0\\ 0 & D^+ \end{array} \right ).
		\]
		[For a nonsingular matrix $M$, $M^+$ denotes the transpose inverse of $M$.]
		Now $A$ and $D$ are similar, hence $\varphi$ and $\psi$ are conjugate.
	\end{proof}
		
\begin{lemma}  \label{lemma-7}
	Let $\varphi \in \OG(V)$ be orthogonally indecomposable of type 1e. Then $\varphi$ is the unipotent factor in the multiplicative Jordan-Chevalley decomposition of a cyclic
	transformation $\eta \in \OG(V,Q)$.	
\end{lemma}

\begin{proof}
	Let $\eta \in \OG(V)$ be cyclic with minimal polynomial $\mu(\eta) = (x^2+x+1)^{2m}$
	(or $(x+\lambda)(x+\lambda^{-1})$ for $K \ne \GF(2)$  and $\lambda \ne 1$). Such transformations exist by Huppert \cite[3.5 Satz]{Huppert-1980b}.
	
	Let $\eta = \eta_{\mathrm{U}} \eta_{\mathrm{S}}$ be the multiplicative
	Jordan-Chevalley decomposition of $\eta$, where 
	$\eta_{\mathrm{U}}$ is unipotent and $\eta_{\mathrm{S}}$ is semisimple.
	Then $\eta_{\mathrm{U}} \in \SpG(V,f_Q)$. By \ref{lemma-5}, $\eta_{\mathrm{U}}$ is orthogonally indecomposable with minimal polynomial $(x+1)^{2m}$ (To see this in the case $\mu(\eta) =(x^2+x+1)^{2m}$, consider an algebraic closure of $K$). 
	Now $Q(v + v\eta_{\mathrm{S}}) = f_Q(v, v +v \eta_{\mathrm{S}})$.
	Hence $\eta_{\mathrm{S}} \in \OG(V)$
	as $\eta_{\mathrm{S}}$ is fixfree. So $\eta_{\mathrm{U}} \in \OG(V,Q)$. By \ref{lemma-6}, $\varphi$ is conjugate to $\eta_{\mathrm{U}}$.
\end{proof}

\section{Proof of theorem \ref{theorem-1}}

Again, we assume that $\ch K = 2$. First we show:
\begin{proposition} \label{prop-3}
	Let $\varphi \in \OG(V)$ be orthogonally indecomposable. Then 
	\begin{enumerate}
		\item $\varphi$ is inverted by an involution $\sigma \in \OG(V,Q)$ with $\dim \Fix(\sigma) = \frac{\dim V}{2}$.
		\item If $\varphi$ is cyclic and $\sigma \in \OG(V,Q)$ is an involution inverting $\varphi$, then $\Fix(\varphi) = \Fix(\sigma) \cap \Fix(        \varphi \sigma)$.
		\item If $\varphi$ is of type 1o or type 2 and unipotent, then 
		$\varphi$ is inverted by  involutions $\rho, \sigma \in \OG(V,Q)$ with $\dim \Fix(\rho) = \frac{\dim V}{2}$ and  $\dim \Fix(\sigma) = \frac{\dim V}{2} + 1$.
	\end{enumerate}
\end{proposition}

\begin{proof}
	First let $\varphi$ be cyclic.  By \cite{EllersNolte-1982}, Lemma 1], $\varphi$ is inverted by an involution $\sigma \in \OG(V)$. Then $\varphi$ is also inverted by the involution $\tau := \varphi \sigma$. We have 
	$\dim \Fix(\sigma), \dim \Fix(\tau) \ge \frac{\dim V}{2}$ as $\ch K = 2$. On the other hand, 
$\Fix(\sigma) \cap \Fix(\tau) \le \Fix(\varphi)$ has at most dimension one. If 
$\varphi$ is not unipotent, then $\dim \Fix(\sigma) = \dim \Fix(\tau) = \frac{\dim V}{2}$. If $\varphi$ is  unipotent, then $\varphi \not \in \SOG(V)$
 so that $\dim \Fix(\sigma) = \dim \Fix(\tau) \pm 1$.
 
 Now let $\varphi$ be of type 1e. We can apply \cite{EllersNolte-1982}, Lemma 4] or use \ref{prop-2}: \\
 Let $\psi \in \OG(V)$ be cyclic with minimal polynomial 
 $(x^2+x+1)^{2t}$, where $\dim V = 4t$. As just seen, $\psi$ is inverted by an involution $\sigma \in \OG(V)$ with $\dim \Fix(\sigma) = 2t$. Clearly, $\sigma$ inverts the
 unipotent factor in the multiplicative Jordan-Chevalley decomposition of $\psi$.
 
 Finally, let $\varphi$ be of type 1o. By \ref{prop-2},
 $\varphi = \psi^2$, where $\psi$ is cyclic and unipotent. Let
 $\psi = \rho \sigma$, where 
 $\rho \in \SOG(V)$ and $\sigma \in \OG(V) -  \SOG(V)$. Then $\varphi = \rho \rho^{\sigma} = \sigma^{\rho} \sigma$.
\end{proof}

It follows immediately from \ref{prop-3} that $\SOG(V)$ is bireflectional if
$\dim V \equiv 0 \mod 4$.
To complete the proof of \ref{theorem-1}(1), it remains to show that $\SOG(V)$ is not bireflectional if
$\dim V = 2m \equiv 2 \mod 4$ and $V$ is not a hyperbolic plane over $\GF(2)$. In this case, $V$ has a basis of anisotropic vectors. It follows from  \ref{remark-3} that there exists a  product $\psi$ of reflections $\sigma_1, \dots, \sigma_{2m} \in \OG(V)$ such that $\Bahn(\psi) = V$. 
 Then $\Fix(\psi) = 0$ so that $\psi$ cannot be bireflectional by \ref{prop-3}.

For the proof of
\ref{old-theorem}(2) in \cite{KN-1987a}, we showed that if $\ch K \ne 2$, then every $\varphi \in \OG(V)$, that is inverted by an $\sigma \in \OG(V)$ admits an orthogonal decomposition $V = U_1 \perp \dots \perp U_t$, where the spaces $U_j$ are $\langle \varphi, \sigma \rangle$-invariant and orthogonally indecomposable.
Unfortunately this result does not remain valid for fields of characteristic 2. Instead, we use a centralizer argument. First we show:

\begin{lemma}                                   \label{lemma-8}
	Let $\varphi \in \OG(V)$.
	Let $T$ be a $\varphi$-invariant totally isotropic subspace of $V$.
	If $\dim T = \frac{1}{2} \dim V$, then $\varphi \in \SOG(V)$.
\end{lemma}

\begin{proof}
	Assume that $m := \dim T = \frac{1}{2} \dim V$. Let $S$ be a totally isotropic complement of $T$. Extend a basis of $T$ by a suitable basis of $S$ to a basis of $V$ such that in this basis
	\[
	\varphi  = \left (\begin{array} {cc} P & X\\ 0 & P^+ \end{array} \right ),
	f_Q = \left (\begin{array} {cc} 0 & \Idm_m\\ \Idm_m & 0 \end{array} \right ),
	\]
	Then $\varphi = \psi \zeta$, where 
	\[
	\psi = \left (\begin{array} {cc} P & 0\\ 0 & P^+ \end{array} \right ),
	\zeta = \left (\begin{array} {cc} \Idm_m & P^{-1}X\\ 0 & \Idm_m \end{array} \right ).  
	\]
	 Clearly, $\det \varphi = 1$. So let $\ch K = 2$. We have to show that $\dim \Bahn(\varphi)$ is even. First we see that $\psi \in \OG(V)$ as  $\psi \in \SpG(V, f_Q)$, $S \psi = S $, and $T \psi = T$. Hence  $\zeta \in \OG(V,Q)$. Now  $\zeta$ is involutory and $\Bahn(\zeta)$ is totally isotropic. It follows from \ref{lemma-1} that every orthogonally indecomposable orthogonal summand of $\zeta$ is of type 1 so that $\dim \Bahn(\zeta)$ is even. Obviously, $\dim \Bahn(\psi)$ is even.
\end{proof}

\begin{lemma}                                   \label{lemma-9}
	Let $\varphi \in \OG(V)$ be unipotent. Assume that all orthogonally indecomposable orthogonal summands of $\varphi$ are of  type $1e$. Then
	\begin{enumerate}
		\item $\Cent(\varphi) \le \SOG(V)$.
		\item If $\alpha \in \OG(V)$ inverts $\varphi$, then
		$\alpha \in \SOG(V)$.
	\end{enumerate}
\end{lemma}

\begin{proof}
	1: Let $T = \rad \Bahn(\varphi) + \rad \Bahn^2(\varphi) + \dots + \rad \Bahn^{\infty}(\varphi)$. Then $T$
	is a totally isotropic, $\Cent(\varphi)$-invariant subspace of $V$ and $\dim T = \frac{1}{2} \dim V$.
	By \ref{lemma-8}, $\Cent(\varphi) \subseteq \SOG(V)$. 
	
	2: By \ref{prop-3}, $\varphi$ is inverted by an involution
	$\sigma \in \SOG(V)$. Then $\alpha \sigma \in \Cent(\varphi) \le \SOG(V)$ so that $\alpha \in \SOG(V)$.
\end{proof}

\begin{lemma}
	Let $\ch K = 2$, $\dim V \equiv 2 \mod 4$.
	Let $\varphi \in \SOG(V)$ be  bireflectional. Then  $\varphi$ has an orthogonal summand of type 1o or a cyclic unipotent  orthogonal summand.
\end{lemma}

\begin{proof}
Let $\sigma \in \OG(V)$ be an involution inverting $\varphi$.
Clearly, $F := \Fix^{\infty}(\varphi)$ and $B := \Bahn^{\infty}(\varphi)$ are $\sigma$-invariant.	Assume  that all unipotent orthogonally indecomposable orthogonal summands of $\varphi$ are of type 1e. Then by \ref{prop-3} and \ref{lemma-11}, $\dim \Fix(\sigma)$ is odd.
\end{proof}
This finishes the proof of theorem \ref{theorem-1}(2).

\section{Proof of  theorem \ref{theorem-2}}

\begin{lemma}                                                                 \label{lemma-10}
	Let $\ch K \ne 2$. If $\varphi \in \OG(V)$ has no orthogonal summand of odd dimension, then 
	$\Cent(\varphi) \subseteq \SOG(V)$.
\end{lemma}

\begin{proof}
	By contraposition.
	Let $\xi \in \Cent(\varphi) - \SOG(V)$. Clearly, $\Fix^{\infty}(-\xi)$ is $\varphi$-invariant. By \ref{prop-1}, 
	$\dim \Fix^{\infty}(-\xi)$ must be  odd. 
\end{proof} 

\begin{lemma}                                   \label{lemma-11}
	Let $\ch K = 2$.
	Let $\varphi \in \OG(V)$ with $\Fix(\varphi) = 0$. Let
	$\xi \in \Cent(\varphi)$.
	Then every unipotent orthogonally indecomposable  orthogonal summand of $\xi$ is of type 1 and $\Cent(\varphi) \subseteq \SOG(V)$.
\end{lemma}

\begin{proof}
	Suppose that there exists a nondegenerate $\xi$-cyclic subspace $U \le \Fix^{\infty}(\xi)$ with
	$\dim U = 2t$. Put
	$T := \Fix^t(\xi) \cap \Bahn^t(\xi)$.
	Then $T \le T^{\perp}$ but by \ref{lemma-1}, $T$ is not totally isotropic.
	On the other hand, $T$ is $\varphi$-invariant. Hence $T = T(\varphi -1)$ so that $T$ must be totally isotropic,
	a contradiction.
\end{proof}

\begin{corollary}                                                                 \label{cor-5}
	Let $\ch K = 2$.
	Let $\varphi \in \OG(V)$. If every unipotent orthogonally indecomposable  orthogonal summand of $\varphi$ is of type 1e, then $\Cent(\varphi) \subseteq \SOG(V)$.
\end{corollary}

\begin{proof}
	Clearly, $\Fix^{\infty}(\varphi)$ and $\Bahn^{\infty}(\varphi)$ are $\Cent(\varphi)$-invariant.
	Apply \ref{lemma-9} and \ref{lemma-11}.
\end{proof}

We finish the proof of \ref{theorem-2}:
Suppose that $\varphi$ is reversible but not bireflectional. Then $\Cent(\varphi) \not \le  \SOG(V)$. But then by \ref{cor-5}, $\varphi$ has a unipotent orthogonally indecomposable summand of type 2 or 1o, in contradiction to theorem \ref{theorem-1}.


\end{document}

\typeout{get arXiv to do 4 passes: Label(s) may have changed. Rerun}